\documentclass[mathpazo]{cicp}
\usepackage[T1]{fontenc}
\usepackage[latin9]{inputenc}
\usepackage[english]{babel}
\usepackage{float}
\usepackage{booktabs}
\usepackage{amstext}
\usepackage{graphicx}
\usepackage{xcolor}

\usepackage{comment}

\makeatletter

\usepackage{babel}
\usepackage{amsmath}
\usepackage{amssymb}

\makeatother

\begin{document}

\title{A wavelet-adaptive method for multiscale simulation of turbulent flows in flying insects}

\author[Engels et~al.]{
Thomas Engels\affil{1,4}\comma\corrauth,
Kai Schneider\affil{2},
Julius Reiss\affil{3},
Marie Farge\affil{4}
}

\address{ \affilnum{1} Institute of Biosciences, University of Rostock, Rostock, Germany \\
          \affilnum{2} Aix-Marseille Universit\'e, CNRS, I2M UMR 7373, Marseille, France \\
          \affilnum{3} ISTA, Technische Universit\"at Berlin, Berlin, Germany \\
          \affilnum{4} LMD UMR 8539 \'Ecole Normale Sup\'erieure-PSL, Paris, France}

\emails{{\tt thomas.engels@ens.fr} (T.~Engels)}

\begin{abstract}
We present a wavelet-based adaptive method for computing 3D multiscale flows in complex, time-dependent geometries, implemented on massively parallel computers. While our focus is on simulations of flapping insects, it can be used for other flow problems. We model the incompressible fluid with an artificial compressibility approach in order to avoid solving elliptical problems. No-slip and in/outflow boundary conditions are imposed using volume penalization. The governing equations are discretized on a locally uniform Cartesian grid with centered finite differences, and integrated in time with a Runge--Kutta scheme, both of 4th order. The domain is partitioned into cubic blocks with different resolution and, for each block, biorthogonal interpolating wavelets are used as refinement indicators and prediction operators. Thresholding the wavelet coefficients allows to generate dynamically evolving grids, and an adaption strategy tracks the solution in both space and scale. Blocks are distributed among MPI processes and the grid topology  is encoded using a tree-like data structure. Analyzing the different physical and numerical parameters allows us to balance their errors and thus ensures optimal convergence while minimizing computational effort. Different validation tests score accuracy and performance of our new open source code, \texttt{WABBIT}. Flow simulations of flapping insects demonstrate its applicability to complex, bio-inspired problems. 
\end{abstract}

\ams{65M06; 65M50; 65M55; 65M85; 65T60; 76M20; 76D05; 76Z10}
\keywords{adaptive numerical method, wavelets, volume penalization, artificial compressibility, flapping flight}

\maketitle

\section{Introduction}
Computing multiscale flows in complex geometries, which may move or deform, is required in numerous 
applications, {\it e.g.}, many biological flow problems such as flying insects or beating hearts. 
This remains a major challenge for computational fluid dynamics, especially in the
turbulent flow regime. To simulate turbulent flows in the presence of moving
boundaries, numerical techniques are needed that allow the solution to be tracked in
both space and scale, and the numerical grid to be dynamically adapted accordingly. If
the fluid-structure interaction must also be taken into account, it is all the more difficult
since the motion of the boundary is no longer known \emph{a priori}, but depends on its
nonlinear interaction with the fluid.
%
In this regard adaptive numerical discretization methods, which can be traced back to the 1980s \cite{Brandt1977,Berger1984}, are indeed attractive. In many cases they can be much more competitive than schemes on uniform fine grids, depending on the character of the solution. However, for adaptive discretizations two major challenges can be identified: their actual implementation on massively parallel supercomputers and the numerical error analysis of adaptivity.

The implementation of the code is crucial to optimize computing, and two conceptually different approaches can be distinguished: one uses point-based techniques, while the other uses block-based techniques. In the former, the error indicator determines for each grid point whether it is significant or not ({\it e.g.}, \cite{Roussel2010,Kevlahan2005}), while in the latter significant grid points are clustered in patches with the drawback of including non-significant points and thus decreasing the compression rate \cite{Domingues2003}. Due to the hardware layout of modern CPU, block-based implementations are in many cases more competitive and have become increasingly popular during the last years, see, {\it e.g.},~\cite{Schornbaum2018} for a recent review on available software packages.  The locally regular block data can be transferred in one contiguous chunk to the CPU cache, which greatly increases the performance despite an increase in the number of points.

The mathematical support for adaptivity needs to provide reliable error estimators of the solution and, for evolutionary problems, a prediction of the grid used to compute the next time step. For both, a variety of heuristic criteria exists, {\it e.g.}, gradient-based approaches~\cite{Deiterding2016}.  Adaptive mesh refinement algorithms use these heuristic criteria and are \emph{error-indicated} methods because of their heuristic nature. \emph{A posteriori} error estimators \cite{BeckerRannacher2001} are mathematically rigorous but require solving expensive adjoint problems. 

Wavelets and related multiresolution analysis techniques provide likewise a mathematical framework and yield reliable error estimators, coupled with high computational efficiency; thus they are well suited for developing adaptive solvers with \emph{error control}.  The idea of wavelet analysis is to decompose data into contributions in both space and scale (and possibly direction). The wavelet transform has been introduced by Grossmann and Morlet~\cite{Grossmann1984}, and the algorithm of the fast wavelet transform by Mallat~\cite{Mallat1989}. Nonlinear approximation~\cite{DeVore1998} provides the conceptual support for adaptivity; indeed, it introduces a systematic way to classify functions according to the sparsity of their representation in wavelet space, {\it i.e.}, the possibility of describing a function by a small number of wavelet coefficients. Wavelet-based algorithms for solving PDEs, starting with~\cite{Liandrat1990}, were the first ones to be proved to guarantee the best $n$-term approximation \cite{Cohen1998}. The books of Cohen~\cite{Cohen2000} and M\"uller~\cite{Mueller2003} and the review articles~\cite{Domingues2011,Schneider2010} give a detailed overview on the subject.

To take into account boundary conditions for complex geometries, in particular at solid walls, the family of immersed boundary methods (IBM)~\cite{Peskin2002,Mittal2005,Schneider2015} is now often preferred over the traditional body-fitted grids for complex geometries. While the latter can yield better accuracy near the wall, IBM approaches convince with their ease of implementation and great flexibility. An early application for computing flow in a beating heart can be found in~\cite{Peskin1977}, illustrating that this class of methods has been designed from the start for complex, bio-inspired flow problems. Among the variety of methods, the volume penalization method~\cite{Arquis1984}, physically motivated by flow in porous media, is furthermore distinguished by its convergence proof for nonlinear Navier--Stokes~\cite{Angot1999}. Moreover, its convergence properties are improved in the turbulent regime~\cite{Nguyenvanyen2014}.

%

In biological applications the flow can be typically described as incompressible, which requires specialized numerical techniques to decouple the pressure from the velocity, most of which involve solving a Poisson equation. While there are adaptive solvers following this strategy, {\it e.g.}, \cite{Kevlahan2005}, alternative methods such as the Lattice--Boltzmann method are an attractive choice for adaptive simulations due to their explicit character~\cite{Schornbaum2018}. Another recently revisited method~\cite{Ohwada2010} is the artificial compressibility method~\cite{chorin1967}, which introduces a large but finite speed of sound and avoids solving elliptic equations.

In this work, we present a wavelet-based fully adaptive approach for solving the incompressible Navier--Stokes equations in complex and time-varying geometries on massively parallel computers. To this end, we combine two physically motivated models: the artificial compressibility method, which relaxes the kinematic incompressibility constraint and allows for a finite speed of sound, and the volume penalization, which models solid objects as porous media with finite permeability. In the limit when the speed of sound and the permeability parameters tend to infinity and to zero, respectively, we recover the incompressible Navier--Stokes equations with no-slip boundary conditions.  The numerical method for solving the resulting model equations is based on centered finite differences and explicit time integration, both of 4th order. The MPI parallel implementation uses block-based grids, which are coarsened and refined using biorthogonal wavelets. While the framework we present is quite general and can be applied to a large variety of fields, our main motivation for its development comes from the spectacular flight capabilities of flying insects, which use flapping wings to produce the required aerodynamic forces. 

The remainder of the paper presents the governing equations (sec.~\ref{sec:govequ}), followed by the numerical method including discretization errors (sec.~\ref{sec:nummeth}), and the multiresolution analysis for introducing adaptivity (sec.~\ref{sec:Adaptivity-and-thresholding}). The choice of the parameters for balancing the different errors is discussed in sec.~\ref{sec:Error-Balancing}. The parallel implementation is presented in sec.~\ref{sec:Parallel-implementation}. Different numerical results are shown in sec.~\ref{sec:results}, and the performance of the code is analyzed in sec.~\ref{sec:Performance}. Finally, an outlook is given in sec.~\ref{sec:outlook} conclusions are drawn in sec.~\ref{sec:conclusions}.

\section{Governing equations and numerical method\label{sec:govequ}}

\subsection{Artificial compressibility}
The applications we have in mind involve exterior, incompressible flows past obstacles of complex shape which may vary in time. Numerical methods for such flows usually require the solution of a Poisson-type, elliptic PDE. Even though powerful libraries ({\it e.g.}, PETSc \cite{abhyankar2018}) exist for that purpose, the solution of the Poisson equation is causing major computational cost, especially on irregular, time-dependent grids.

Instead of the incompressible Navier--Stokes equations,
\begin{eqnarray}
\partial_{t}\underline{u}  +\left(\underline{u}\cdot\nabla\right)\underline{u}+\nabla p-\nu\nabla^{2}\underline{u} & = & 0\label{eq:INS-1}\\
\nabla\cdot\underline{u} & = & 0 \label{eq:INS-2}
\end{eqnarray}
for the velocity $\underline{u}$ and the pressure $p$ in the fluid domain $\Omega_f \subset \mathbb{R}^D$ for $t>0$ and $D=2$ or $3$, completed with suitable boundary and initial conditions, and where $\nu$ is the kinematic viscosity, we consider the artificial compressibility method (ACM). An artificial speed of sound $C_{0}\gg |\underline{u}|$ is introduced and eqn.~(\ref{eq:INS-2}) is replaced by
\begin{equation}
\partial_{t}p + C_{0}^{2}\nabla\cdot\underline{u} +C_{\gamma}p = 0\label{eq:ACM-2}\\
\end{equation}
where $C_{\gamma} >0$ is a dashpot damping term discussed below. Note that the fluid density $\varrho_{f}$ is normalized to unity and the equations are written in dimensionless form. The basic idea of the ACM can be traced back to Chorin~\cite{chorin1967}. More recently, several authors have reconsidered it. Guermond and Minev state that ``the artificial compressibility method may not have been given all the attention it deserves in the literature'' \cite{Guermond2015}. The spirit adopted here, which is to explicitly solve eqns. (\ref{eq:INS-1},\ref{eq:ACM-2}), is inspired by the work of Ohwada and Asinari~\cite{Ohwada2010}. The ACM solution converges to the solution of the incompressible Navier--Stokes equations (INC) with
\begin{equation}
e_{\mathrm{ACM}}=\left\Vert \underline{U}_{\mathrm{ACM}}-\underline{U}_{\mathrm{INC}}\right\Vert _{2}=\mathcal{O}(C_{0}^{-2})\label{eq:ACM-model-error}
\end{equation}
as $C_{0}\rightarrow\infty$. Here, $\underline{U}$ is the state vector, containing $\underline{u}$ and $p$. We will confirm this result numerically later. Note that \cite{Ohwada2010} proposes a strategy to increase the accuracy to $\mathcal{O}(C_{0}^{-4})$ using Richardson extrapolation, but we do not follow this approach because the second order is sufficient when combined with the volume penalization method. We also point out that there is a natural upper limit for $C_{0}$, which is the true speed of sound of the fluid ({\it e.g}. $C_{0}=340\,\mathrm{m}\mathrm{s}^{-1}$ in air), and that any fluid is at least slightly compressible.

\subsection{Volume penalization}

The no-slip boundary conditions at the fluid--solid interface are imposed using the volume penalization method~\cite{Arquis1984,Angot1999}, which is based on the physical intuition of approximating a solid obstacle as
a porous medium with small permeability. The penalized versions of eqns. (\ref{eq:INS-1},\ref{eq:ACM-2}) are solved in the larger domain $\Omega = \Omega_f \bigcup \Omega_s$, where $\Omega_s$ is the solid domain, and read
\begin{eqnarray}
\partial_{t}\underline{u} +\underline{u}\cdot\nabla\underline{u} 
+\nabla p 
-\nu\nabla^{2}\underline{u}
+\frac{\chi}{C_{\eta}}\left(\underline{u}-\underline{u}_{s}\right)
+\frac{\chi_{\mathrm{sp}}}{C_{\mathrm{sp}}}\left(\underline{u}-\underline{u}_{\infty}\right) & = & 0 \label{eq:u_eqn-1}\\
\partial_{t}p +C_{0}^{2}\nabla\cdot\underline{u}+C_{\gamma}p+\frac{\chi_{\mathrm{sp}}}{C_{\mathrm{sp}}}\left(p-p_{\infty}\right) & = & 0,\label{eq:p_eqn-1}
\end{eqnarray}
Here $\chi$ is the indicator function, $C_{\eta}\ll1$ is the penalization
parameter (permeability) and $\underline{u}_{s}$ the solid body velocity
field. The indicator function is $\chi=1$ inside $\Omega_s$ and $\chi=0$
in $\Omega_f$, and we use a thin smoothing layer with thickness proportional
to the grid spacing in the case of moving~\cite{Kolomenskiy2009}
or deforming obstacles~\cite{Engels2012a}. For details on the generation of $\chi$, see~\cite{Engels2015a}. From the penalization
term aerodynamic forces are obtained by volume integration \cite{Angot1999},
\[
\underline{F}=C_{\eta}^{-1}\int_{\Omega}\chi\left(\underline{u}-\underline{u}_{s}\right)\mathrm{d}V,
\]
for which we use the midpoint quadrature rule. We anticipate that the penalization error converges as 
\begin{equation}
e_{\mathrm{VPM}}=\left\Vert \underline{U}_{\mathrm{ACM},{C_{\eta}}}-\underline{U}_{\mathrm{ACM}}\right\Vert _{2}=\mathcal{O}(C_{\eta}^{1/2})\label{eq:penal_error}
\end{equation}
for $C_{\eta}\rightarrow0$ using similar arguments as for INC.

For our application to unbounded exterior flows, an additional sponge
penalization term (index $\mathrm{sp}$) has been added \cite{Engels2012a,Engels2014},
which imposes mean flow conditions, $\underline{u}_{\infty}$ and $p_{\infty}$,
in a layer at the outer boundary of the computational domain. Contrarily
to the penalization for the obstacle itself, this term also acts on
the pressure $p$. We can understand the ACM as a transport of divergence,
which is created both at the boundary and in the nonlinear term. Hence,
some outflow conditions are required to remove it from the physical domain, and
the sponge term here plays the dual role of absorbing outgoing pressure
waves as well as damping outgoing vortical structures. 

In the case of impulsively started flow, the initial condition is
incompatible with the boundary conditions, which causes a singularity
in the aerodynamic force. A strong pressure wave of width $\delta_{s}=\mathcal{O}(C_{\eta}C_{0})$
is created at the boundary, because the penalization requires a time $\mathcal{O}(C_{\eta})$
 to damp the velocity inside the obstacle,
and the ACM transports it with a speed $C_{0}$. Its intensity decays
with the distance $R$ as $R^{-(D-1)}$, where $D=2$ or $3$ is the dimensionality of the
problem. Even though, physically speaking, the pressure wave is an
artifact of the ACM. It still needs to be resolved, as it transports information needed to create an (approximative) solenoidal flow field in the full domain. If $\delta_{s}$
is too small, dispersion errors occur, and spurious oscillations in
the forces can appear. In practice, we ensure that $\Delta x \leq \delta_{s}$.
In order to minimize reflections, the width of the sponge layer $L_{sp}$
is set larger than characteristic length scales of the geometry; a
wider sponge reduces reflections. We then adjust the sponge constant
$C_{sp}$ such that the time the pressure wave spends in the sponge is
equivalent to $\tau$ times the relaxation time, yielding $C_{sp}=L_{sp}(C_{0}\tau)^{-1}$.
Here we set $\tau=20$, but found that the specific choice of both $\tau$
and $L_{sp}$ only has a weak impact on the solution (see Appendix D in suppl. material for details).  

In simulations with periodic boundary conditions, the sponge technique cannot
be used. For such cases, we use the dashpot damping term $C_{\gamma}p$
proposed in \cite{Ohwada2010} to mitigate pressure waves; the constant
is then $C_{\gamma}=1$, as suggested in \cite{Ohwada2010}.

The combination of artificial compressibility and volume penalization
constitutes our physical model. It has four parameters ($C_{\eta}$,
$C_{0}$, $C_{\mathrm{sp}}$, $C_{\gamma}$), the choice of which
will be discussed in detail below. 

\subsection{Numerical method and discretization errors\label{sec:nummeth}}

The model equations (\ref{eq:u_eqn-1}-\ref{eq:p_eqn-1}),
are discretized in space and time. Spatial derivatives
are approximated using centered higher-order finite differences on a collocated, periodic grid.
We use the optimized fourth order scheme proposed by Tam and Webb~\cite{Tam1993} for first
order derivatives (Appendix~A in suppl. material) and a classical, centered, fourth order stencil for second order derivatives. 

Time integration is done with a classical Runge--Kutta-4 scheme.
The explicit scheme implies stability restrictions on the time step
$\Delta t$, namely 
\[
\Delta t<\min\left(\mathrm{CFL}\Delta x_{\mathrm{min}}/u_{\mathrm{eig}},\:\mathrm{CFL}_{\eta}C_{\eta},\:\text{\ensuremath{\mathrm{CFL}_{\nu}}}\Delta x_{\mathrm{min}}^{2}/\nu\right),
\]
where the characteristic velocity is $u_{\mathrm{eig}}=\left|u_{\mathrm{max}}\right|+\sqrt{u_{\mathrm{max}}^{2}+C_{0}^{2}}$,
which we can approximate by $u_{\mathrm{eig}}\approx C_{0}$ as $C_{0}\gg |u_{\mathrm{max}}|$.
The constants are set to $\mathrm{CFL}=\mathrm{CFL}_{\eta}=1.0$ and
$\mathrm{CFL}_{\nu}=1/4$. We therefore expect the fully discretized
solution to converge with 4th order accuracy in space and time,
\begin{align}
e_{\mathrm{FDM}}=\left\Vert \underline{U}_{\mathrm{ACM}}^{\Delta x, \Delta t}-\underline{U}_{\mathrm{ACM}}\right\Vert _{2} & =\mathcal{O}(\Delta x^{4}) + \mathcal{O}(\Delta t^{4}),\label{eq:discretization_error_4th}
\end{align}
which we will verify later numerically. This is valid in the case without penalization. As pointed out in \cite{Nguyenvanyen2014}, the discretization
(eq. \ref{eq:discretization_error_4th}) and penalization (eq. \ref{eq:penal_error})
error cannot be treated independently. A too small value for $C_{\eta}$
results in a loss of regularity of the exact solution and by consequence the discretization error increases, thus there is  an optimal choice of $C_{\eta}$ as a function of $\Delta x$. As discussed in our previous work \cite{Engels2014,Engels2015a},
the relation is
\begin{equation}
C_{\eta}=\left(K_{\eta}^{2}/\nu\right)\Delta x^{2}\label{eq:ceta-scaling}
\end{equation}
where $K_{\eta}$ is a constant which depends only on the discretization
scheme. It can be interpreted as the number of grid points in the
penetration boundary layer inside the obstacle. Using eqn.~(\ref{eq:ceta-scaling}),
first order convergence \cite{Nguyenvanyen2014} is assured for any
value of $K_{\eta}$, {\it i.e.},
\begin{equation}
\left\Vert \underline{U}_{\mathrm{ACM},C_\eta}^{\Delta x, \Delta t}-\underline{U}_{\mathrm{ACM}}\right\Vert _{2}=\mathcal{O}(\Delta x)+\mathcal{O}(\Delta t^{4}),\label{eq:VPM_error}
\end{equation}
as $\Delta x,C_{\eta}\rightarrow0$. The error offset can be tuned
by choosing $K_{\eta}$. In the present work, we obtained best results
using the same values of $0.05\leq K_{\eta}\leq0.4$ as reported in
\cite{Engels2015a} with a Fourier discretization, which suggests
a certain universality of this scaling.

\section{Multiresolution analysis and grid adaptation\label{sec:Adaptivity-and-thresholding}}

Motivated by the multitude of scales which are excited in the flow, especially in the turbulent regime, and close to the boundary, whose location may even change in time, we introduce dynamically adaptive grids. The grid is supposed to track automatically the motion in both scale and position and allows furthermore its prediction for the next time step.  To this end we use multiresolution analysis and decompose the solution into biorthogonal wavelets.
First we recall the point-value multiresolution analysis and we give the link with biorthogonal wavelets. Some drawbacks of point-value multiresolution analysis are pointed out and we show how they can be overcome using lifted wavelets. Some details on the parallel implementation are given later in section~\ref{sec:Parallel-implementation}.

For ease of presentation we limit the description to the one-dimensional scalar-valued case, denoted by $u$. The extension to higher dimension and vector-valued data is outlined at the end of the section.

\subsection{Point-value multiresolution analysis}

Discrete multiresolution analysis introduced by Harten~\cite{Harten1993,Harten1995,Harten1996} is well-suited  for introducing adaptivity into discretization schemes of PDEs. In particular for finite difference methods the point-value multiresolution~\cite{Harten1993}, which is directly related to the interpolating subdivision scheme by Deslauriers \& Dubuc~\cite{Deslauriers1987,Deslauriers1989}, has been designed.

Starting point is a nested sequence of uniform grids at level $J$, $X^{J}=\{x_{i}^{J}\}_{i=0}^{2^{J}}$ with an odd number of grid points defined by $x_{i}^{J}=i/2^{J}$ for $i=0,\ldots2^{J}$.
The grids satisfy $x_{2i}^{J+1}=x_{i}^{J}$ for $i=0,...,2^{J}$ and $x_{2i-1}^{J+1}=(x_{i}^{J}+x_{i-1}^{J})/2$ for $i=1,...,2^{J}$.

For the solution $u$ given on the grid $X^J$ we can then define a prediction operator $P_{J\rightarrow J+1}:u(X^{J})\rightarrow u(X^{J+1})$
which interpolates values of $u$ onto the next finer grid. 
For this, we use the fourth order interpolating subdividison scheme \cite{Deslauriers1987,Deslauriers1989}
to interpolate missing values at $x_{2i+1}^{J+1}$, {\it i.e.},
\begin{equation}
u_{2i+1}^{J+1}=-\frac{1}{16}u_{i-1}^{J}+\frac{9}{16}u_{i}^{J}+\frac{9}{16}u_{i+1}^{J}-\frac{1}{16}u_{i+2}^{J}.
\end{equation}

Following Harten, we also define a restriction operator $R_{J+1\rightarrow J}:u(X^{J+1})\rightarrow u(X^{J})$
which is simply the decimation by a factor of two, {\it i.e.}, $u_{i}^{J}=u_{2i}^{J+1}$.
Note that $R_{J+1\rightarrow J}\circ P_{J\rightarrow J+1}=\mathrm{Id}$.
Combining restriction and prediction, detail coefficients can be computed by subtracting the values predicted from the coarse grid from the values at the fine grid. At even grid points the details are zero, while at odd grid points we obtain
%
%
\begin{equation}
\widetilde u_{2i+1}^{J}=u_{2i+1}^{J+1}-P_{J\rightarrow J+1}R_{J+1\rightarrow J}u_{2i+1}^{J+1}\label{eq:details}
\end{equation}
%
%
Thus, the solution on the fine grid $u^{J+1}$ can be represented as a solution at the coarser grid $u^{J}$ plus detail coefficients $\widetilde{u}^{J}$. Iterating this procedure from $J$ down to $1$ yields the multiresolution transform. The fine grid solution can be reconstructed by inverting the above procedure. Detail coefficients are mostly significant in regions of steep gradients and discontinuities, and are small or even do vanish in regions where the solution is smooth. They ``measure the local deviation from a polynomial'' and are obtained as the interpolation error.
Thus adaptivity can be introduced by removing small detail coefficients without loosing the precision of the computation.

\subsection{Biorthogonal wavelets}

The point-value multiresolution analysis is intimately related to biorthogonal wavelets, as already pointed out by Harten~\cite{Harten1993},
and the solution $u \in C^s(\mathbb{R})$ ($s \ge 0$) can be expanded into a biorthogonal wavelet series,
\begin{equation}
u(x)=\sum_{i\in\mathbb{Z}}\left\langle u,\widetilde{\phi_{i}^{0}}(x)\right\rangle \phi_{i}^{0}(x)+\sum_{J=0}^{\infty}\sum_{i\in\mathbb{Z}}\left\langle u,\widetilde{\psi_{i}^{J}}(x)\right\rangle \psi_{i}^{J}(x)
\label{eqn:biorthwlseries}
\end{equation}
Note that for a given primary scaling function $\phi$, different dual scaling functions $\widetilde{\phi}$ can be constructed (and vice versa) and thus the choice is not unique and there is some flexibility compared to the orthogonal case (in terms of symmetry, number of vanishing
moments, filter length, etc.).
For orthogonal wavelets we have $\widetilde{\phi}=\phi$ and $\widetilde{\psi}=\psi$. 

For the point-value multiresolution the scaling function $\phi$ (cf. appendix~C  in suppl. material)
is identical to the Deslauriers--Dubuc fundamental function~\cite{Deslauriers1987,Deslauriers1989}.
It can be obtained by interpolating a delta Dirac impulse (using Lagrange interpolation), and the filters are defined to be symmetrical and have minimal support for a given polynomial order.
They are equivalent to the autocorrelation function of orthogonal Daubechies~\cite{Daubechies1988} scaling functions \cite{Mallat2009}. 
The dual scaling functions are delta Dirac distributions $\widetilde{\phi}_{k}^{J}(x)=\delta(x-x_{k}^{J})$.
The wavelets are interpolating (`interpolets'), yield an interpolating subdivision scheme, and correspond to a shifted scaling function $\psi(x)=\phi(2x-1)$, while the dual wavelets are linear combinations of $\delta$ distributions. The construction of these interpolating wavelets was also proposed by Donoho~\cite{Donoho1992} and rediscussed later by Sweldens~\cite{Sweldens1997}.
Note that the primary scaling function is normalized with respect to the $L^{\infty}$ norm, {\it i.e.}, $\phi_{i}^{J}(x) = \phi(2^J x -i)$ implying $\left\Vert \phi_{i}^{J} \right\Vert_{\infty} = 1$, and not with respect to the $L^2$ norm, as it is typically the case.

Both scaling functions and wavelets fulfill refinement relations,
\begin{eqnarray}
\phi(x)=\sum_{n\in\mathbb{Z}}h_{n}\phi(2x-n) \quad \quad \widetilde \phi(x)=\sum_{n\in\mathbb{Z}}\widetilde h_{n}\widetilde \phi(2x-n) \\
\psi(x)=\sum_{n\in\mathbb{Z}}g_{n}\phi(2x-n) \quad \quad \widetilde \psi(x)=\sum_{n\in\mathbb{Z}}\widetilde g_{n}\widetilde \phi(2x-n)
\end{eqnarray}
and the filter coefficients are coupled via $\widetilde g_n = (-1)^{1-n} h_{1-n}$ and $g_n = (-1)^{1-n} \widetilde h_{1-n}$.

Rearranging the restriction and prediction operations of the multiresolution transform yields the filter coefficients. Then applying the fast wavelet transform with these filters and downsampling  will yield equivalent results. The detail coefficients correspond to the wavelet coefficients. 
For the point-value multiresolution we have $\widetilde h_0 = 1$ and $\widetilde h_n = 0$ for $n \ne 0$ and thus $g_1= 1$ and $g_n = 0$ for $n \ne 1$.
The coefficients of $h$ and $\widetilde g$ depend on the interpolation order, and in the linear case we have $h_n = \{ 1/2, 1, 1/2\}_{n = -1, 0, 1}$
and $\widetilde g_n = \{ -1/2, 1, -1/2\}_{n = 0, 1, 2}$.

The filter $\widetilde h$ illustrates the loose downsampling of the data (simple decimation by a factor of two) when going to coarser scale and it is thus not a low pass filter in the classical sense as no high frequency contributions are eliminated. 
The corresponding scaling function $\widetilde \phi$ is a distribution, and even not in $L^2(\mathbb{R})$.
The wavelet $\psi$ is a shifted scaling function and hence not a wavelet with vanishing mean.
The corresponding filter $g$ is thus not a band pass filter neither. 

These wavelets resulting from Harten's point-value multiresolution are Cohen--Daubechies--Feauveau
(CDF) biorthogonal wavelets~\cite{Cohen1992} which are symmetric and compactly supported. For cubic interpolation we have CDF 4/0 wavelets, because the analyzing wavelet has $4$ vanishing moments while the synthesizing one has zero vanishing moments. The corresponding filter coefficients $h,g$ and $\widetilde{h},\widetilde{g}$ are summarized in appendix C (suppl. material) together with an illustration of the functions. 

This biorthogonal wavelet decomposition yields good compression properties of the data, as the analyzing wavelet $\widetilde \psi$ has vanishing moments. However, there is no separation into scales with well separated frequency bands, as the loosely downsampling does not filter out high frequency contributions (which produces aliasing) and, since the synthesizing wavelet has no vanishing moment, all frequencies are mixing up. 
For denoising~\cite{Donoho1994} and Coherent Vorticity Simulation~\cite{Farge1999} such decompositions are only of limited use and for this reason we are applying Sweldens' lifting scheme to construct second generation wavelets. The CDF 4/0 biorthogonal wavelets are upgraded to CDF 4/4~\cite{Sweldens1997,Sweldens1998,Sweldens2000} which yield a reasonable scale and frequency separation.

\subsubsection{Lifted biorthogonal wavelets}
To improve the frequency selectivity of the point-value multiresolution decomposition we apply (primary) lifting to the wavelet $\psi$ and thus construct a wavelet with vanishing moments. To this end we modify the primary wavelet $\psi$ and the dual scaling function $\widetilde \phi$ correspondingly, while the secondary wavelet $\widetilde \psi$ and the primary scaling function $\phi$ remain unchanged. For the filters the lifting operations are applied accordingly~\cite{Sweldens1998}.
More precisely we determine lifting coefficients $s_{n}$ such that
\begin{eqnarray}
\psi(x) \, & =& \, \psi^\mathrm{old}(x) \, - \, \sum_n s_n \phi^\mathrm{old}(2x -n) \\
\widetilde \phi(x) \, & =& \, \sum_n  \widetilde h_n \widetilde \phi^\mathrm{old}(2x -n) \, + \, \sum_n s_n \psi^\mathrm{old}(x -n) 
\end{eqnarray}
requiring that the new $\psi$ has $M$ vanishing moments, {\it i.e.}, $\int x^p \psi(x) dx = 0$ for $p=0,...,M-1$.
In the linear interpolation case, {\it i.e.}, CDF 2/0, we have $s_0 = s_1 = -1/4$ and for the resulting lifted filters we get $\{g_n, n=-3,...,1\} = \{/1/8, -1/4, 3/4, -1/4,\\ -1/8\}$
and $\{ \widetilde h_n, n=-2,...,2\} = \{/1/8, -1/4, -3/4, -1/4, 1/8\}$, yielding thus CDF 2/2 wavelets, which have two vanishing moments.
Applying lifting to the cubic interpolatory wavelets CDF 4/0 yields CDF 4/4. The functions $\phi$ and $\psi$ as well as their dual functions and filter coefficients are illustrated in Appendix~C (suppl. material).


%
%

\subsection{Nonlinear approximation and thresholding}
The biorthogonal wavelet representation of the solution $u$ (eqn.~\ref{eqn:biorthwlseries}) can be used for introducing an adaptive grid.
To this end we apply thresholding to the detail coefficients $\widetilde u^{J}_{i} = \left\langle u,\widetilde{\psi_{i}^{J}}(x)\right\rangle$ and mark details with magnitude below a constant, level-independent threshold $C_{\varepsilon}$, precisely if $|\widetilde u^{J}_{i}| / \left\Vert u \right\Vert_\infty < C_\varepsilon$.
As a detail coefficient $\widetilde u^J_i$ is directly related to a grid point $x^{J+1}_{2i+1}$, we could remove the corresponding grid points. In our block-based structure, a block is coarsened only if \emph{all} its details are smaller than $C_\varepsilon$, {\it i.e.}, some insignificant points can be included.
The thresholding error can be estimated in the $L^\infty$ norm and is directly proportional to the threshold, {\it i.e.}, after $N_t$ time steps we have,
\begin{equation}
\left\Vert u_{\mathrm{ACM}}^{\Delta x,C_\varepsilon}-u_{\mathrm{ACM}}^{\Delta x}\right\Vert _{\infty}\overset{\mathrm{def}}{=}e_{\mathrm{MR}}\leq N_{t}C_{\varepsilon}\label{eq:thresholding-error}
\end{equation}
as shown in~\cite{Cohen2000b,Roussel2003}. 
This upper bound assumes a worst-case scenario that the filtering introduces an error $\mathcal{O}(C_{\varepsilon})$ in all $N_{t}$ time steps, which of course depends on the problem. For example in pure transport, the production of fine scales is related only to the lack of translation invariance of the multiresolution transform. 
Eventually, if no velocity is imposed, the error is in fact $e_{\mathrm{MR}}\leq C_{\varepsilon}$.
Thus, whether the number of time steps enters eqn.~(\ref{eq:thresholding-error})
depends on the non-linearity. Assuming the pure transport-based CFL condition, we have $N_{t}\propto\Delta x^{-1}$.



\subsection{Extension to higher dimensions and vectors}
Extensions in dimension larger than one are obtained by tensor product of the one-dimensional wavelets. For multiresolution constructions in two dimensions we obtain wavelets in the horizontal, vertical and diagonal directions, while in three dimensions we have seven directions. For details we refer to~\cite{Mallat2009}. For vector-valued quantities, {\it e.g.}, the state vector $\underline U$, the biorthogonal wavelet expansion is applied to each component of the vector and thus a vector-valued wavelet series is obtained. Hence the coefficients become vector-valued, but not the basis functions. During thresholding, each component is normalized by the $L_\infty$-norm of the total (global) field.

\subsection{Grid prediction and time stepping algorithm\label{subsec:Wavelet-Based-Multiresolution-Alg}}
The algorithm to advance the numerical solution $\underline{U}\left(t^{n},\underline{x}\right)$ (state vector containing $\underline{u}$ and $p$) on the block-structured grid $\mathcal{G}^{n}$ to the new time level $t^{n+1}$ is as follows: 
\begin{enumerate}
\item Refinement. Given $\underline{U}\left(t^{n},\underline{x}\right)$
on the grid $\mathcal{G}^{n}$, we first extend the grid by adding
neighbors in position and scale. This concept of a safety zone, introduced
in~\cite{Liandrat1990}, ensures that the solution is tracked in scale and space. Adding neighbors in scale captures newly emerging scales due to the quadratic nonlinearity, which can produce scales twice as small in one time step. 
The grid $\mathcal{G}^{n}$ is thus extended to $\mathcal{\widetilde{G}}^{n}$ by refining all blocks by one level.
\item Evolution. On the extended grid $\mathcal{\widetilde{G}}^{n}$, we
solve equations~(\ref{eq:u_eqn-1}--\ref{eq:p_eqn-1}) using finite
differences with explicit time-marching. 
\item Coarsening. We compute the detail coefficients of $\underline{U}(t^{n+1},\underline{x})$
for each block, and tag it for coarsening if the largest detail is
smaller than the prescribed threshold $C_\varepsilon$. 
Next, the tag is removed if coarsening would result in a non-graded grid, and finally the tagged blocks are merged with their sister blocks. This procedure is repeated until no more blocks can be coarsened.
\end{enumerate}
This well established algorithm (see, {\it e.g.}, \cite{Liandrat1990,Schneider2010,Domingues2003,Kevlahan2005,Kolomenskiy2016a})
is at the core of our adaptive strategy. 
In practice, the above algorithm is modified by imposing an upper
bound for the refinement, {\it i.e.}, a level $J_{\mathrm{max}}$ is introduced at which no refinement is allowed anymore. 
The largest possible value in the current implementation is $J_{\mathrm{max}}=18$. 
This is done not only because of finite computational resources, but also in the interest of balancing the discretization- and thresholding error. While the former error depends mainly on $\Delta x$ (and thus directly on $J_{\mathrm{max}}$), the latter is determined by $C_{\varepsilon}$. The relation between both is discussed below.

If the maximum level $J_{\mathrm{max}}$ is reached after the refinement
stage and the detail coefficients are significant after the evolution
step, we can choose between two possibilities: either, the corresponding
block is coarsened, even though its details are significant, or it
remains at $J_{\mathrm{max}}$. The former option corresponds to dealiasing
and avoids accumulation of energy in high wavenumbers, if viscous dissipation is not sufficient to stabilize the simulation. In low Reynolds number simulations, it can be omitted for the benefit of improved precision, while in most simulations it is required for stabilization.

While the volume penalization term does usually result in large detail
coefficients near the fluid-solid interface, it can happen that
for large values of $C_{\varepsilon}$ the grid is coarsened at the
interface. In this case, oscillations appear in the aerodynamic forces,
which are undesirable. Therefore, we use the mask function $\chi$
as secondary criterion for block coarsening, {\it i.e.}, a block is not coarsened if it satisfies $\max(\chi) - \min(\chi) > 0$. Consequently, blocks containing the fluid-solid interface are coarsened only if the dealiasing option is used. The initial grid at $t=0$ is created in an iterative process, starting from the coarsest level $J_\mathrm{min}$ and likewise takes $\chi$ into account.

\subsection{Balancing modeling, discretization and thresholding errors\label{sec:Error-Balancing}}

The numerical solution of the Navier--Stokes equations with artificial compressibility on adaptive grids and volume penalization involves five parameters in total, ACM ($C_0$), penalization ($C_\eta$), time ($\Delta t$), space ($\Delta x$) and thresholding ($C_{\varepsilon}$).
The accuracy of the computation with respect to the reference solution, {\it i.e.}, the exact solution of the incompressible Navier--Stokes equation, can be estimated from the various error contributions, namely, the artificial compressibilty (eq.~\ref{eq:ACM-model-error}), penalization (eqn.~\ref{eq:penal_error}), discretization (eqn.~\ref{eq:discretization_error_4th})
and thresholding (eqn.~\ref{eq:thresholding-error}) errors and then
applying the triangle inequality, we find,
\begin{equation}
\left\Vert U_{\mathrm{ACM}, C_\eta}^{\Delta x, \Delta t, C_\varepsilon}-U_{\mathrm{INC}}\right\Vert _{2}\, \leq \,  
\mathcal{O}(C_{0}^{-2}) + \mathcal{O}(C_{\eta}^{1/2}) + \mathcal{O}(\Delta x^{4}) + \mathcal{O}(\Delta t^{4}) +  \mathcal{O}(C_{\varepsilon})
\label{eq:allerrors}
\end{equation}
%
%
A suitable choice of the parameters must take their respective scaling into account in order to balance the different errors.
In sec.~\ref{sec:nummeth} we also already discussed the balancing of discretization and penalization errors, which formally reduces the convergence in $\Delta x$ to first order, (eqn.~\ref{eq:VPM_error}). 
Requiring all errors to decrease simultaneously by the same amount, we find in the case without penalization the following relations
\begin{align}
C_{\varepsilon} & \propto\Delta x^{5} & C_{0} & \propto\Delta x^{-2},\label{eq:scaling-no-penalization}
\end{align}
while in the penalized case we have
\begin{align}
C_{\varepsilon} & \propto\Delta x^{2} & C_{0} & \propto\Delta x^{-1/2} & C_{\eta} & =\left(K_{\eta}^{2}/\nu\right)\Delta x^{2},\label{eq:scaling-penalized}
\end{align}
where $0.05\leq K_{\eta}\leq0.4$. In addition, as described above,
we fix the sponge constant to $C_{\mathrm{sp}}=L_{\mathrm{sp}}/C_{0}\tau$
with $\tau=20$ (see Appendix D in suppl. material for details).

\section{Parallel implementation\label{sec:Parallel-implementation}}
\subsection{Data structure\label{subsec:Datastructure}}
Previous work on wavelet-based adaptivity exploited the fact that
the fast interpolating wavelet transform can be used to assign one wavelet coefficient to each grid point, and hence allows judging for each point if its wavelet coefficient is important or can be discarded~\cite{Holmstroem1999}. 
The resulting grids are as sparse as possible, and each point is a node in a tree data structure~\cite{Roussel2003} or hash-table~\cite{Mueller2003}.
The difficulty is that these data structures are poorly aligned with
the memory access of typical modern CPUs. 
Moreover it creates administration overhead for each grid point. 
If the CPU requires a data point
for computing, the vicinity of this point in the 1D memory is also
transferred to the CPUs cache, which explains the high performance
of structured grids.
To turn the achievable sparsity into CPU time
compression, a block-based hybrid data structure is a suitable choice.
In the context of wavelet-based multiresolution, block-based grids
have first been proposed by Domingues et al.~\cite{Domingues2003}
and later, independently, in~\cite{Rossinelli2011a}. The basic idea
is however older; for adaptive computations without wavelet-support,
they have been pioneered by Berger and Collela~\cite{Berger1988}
and further developed by Deiterding and coworkers~\cite{Deiterding2003,Deiterding2005,Deiterding2009}.
An excellent overview over recently developed software packages can
be found in~\cite{Schornbaum2018}.

\begin{figure}
\begin{centering}
\includegraphics[width=.8\textwidth]{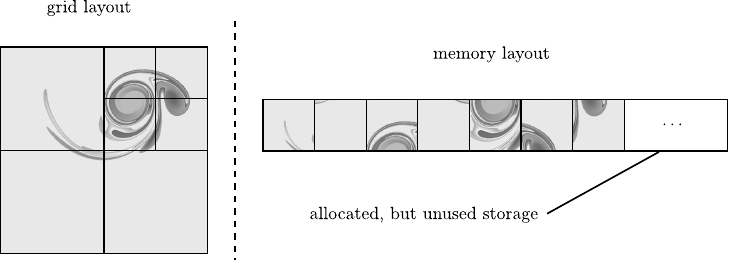}
\par\end{centering}
\caption{Data structure. Left: grid layout, composed of $N_{b}=7$ blocks of the same size but at different positions and levels. Right: memory layout, the grid is stored as simple, contiguous array. Memory allocation done done once on startup, but only the portion required at time $t$  is used. \label{fig:Datastructure}}
\end{figure}

Block-based grids are locally structured and for the relations between blocks, a flexible data structure has to be used. 
In our implementation, we use a tree-based structure
for this task, and define blocks of size $B_{s}^{D}$. All blocks have the
same number of points, which is the same in all directions in the present computations. 
The number of entries in the tree data structure, and thus its overhead, is reduced
by the same factor $B_{s}^{D}$. On the other hand, the total number
of grid points $N_{\mathrm{points}}$ is increased, until for (very) large $B_{s}$, an equidistant
monoblock grid is recovered. Figure~\ref{fig:Datastructure} illustrates
the concept of our data structure. The grid layout contains blocks of
the same size on different levels, while the memory layout is a simple,
contiguous array of size $B_{s}\times B_{s}\times B_{s}\times N_{b}$.
Since the number of blocks $N_{b}$ changes during the simulation,
our code allocates an array with the largest possible number of blocks
for the given memory, which is specified in the program call using
the flag {\tt {-}{-}memory=40.0GB}. Memory allocation is done once
during the initialization and is not deallocated again until the termination of the code. 
This type of memory management is suitable for high-performance
supercomputers, where typically only one application is running on
a given processor.

Since the number of blocks $N_{b}$ in the tree data structure is typically
low ($N_{b}=\mathcal{O}(10^{5})$ blocks, while $N_{\mathrm{points}}=N_{b}B_{s}^{D}=\mathcal{O}(10^{10})$
in 3D applications) there is no need for an external high-performance tree library,
such as \texttt{p4est}~\cite{Burstedde2011}. Here, we design a simple tree library based on the
work of Gargantini~\cite{Gargantini1982,Gargantini1982b}, where octal
numbers are used to encode the position of an element in the tree.
From the octal number,
level and position of the block can be computed, as well as the IDs
of neighboring blocks. It constitutes, along with a refinement flag
and level indicator, the metadata, which is separated from the actual
block data named heavy data. The metadata is stored redundantly on
all CPUs and is kept synchronized, so that neighbor lookup is a local
process, while the heavy data is process-local.

\subsection{Parallelization aspects\label{subsec:Parallelization-Aspects}}
Parallelization appears in our implementation with two functions: to ensure
that each MPI rank has the same amount of blocks to compute (load balancing), and to synchronize the layer of ghost nodes.

The creation and removal of blocks during the simulation can lead
to an imbalanced distribution of blocks among processes. As we use
blocks of identical size, load balancing is guaranteed if the number
of blocks on each process is the same. The load balancing process
first computes the index of all blocks on a space-filling Hilbert-curve~\cite{zumbusch2003,Hilbert1891}, then distributes contiguous chunks
of blocks to the processes. Special care is devoted to keeping the
required number of transfers as small as possible. The Hilbert-curve
preserves the locality of blocks and minimizes the interface between
processes.

Any parallelization by decomposition of the domain requires communication
between the processes, either in the form of flux transfer or ghost
node synchronization. Our blocks are equipped with a layer of $n_{g}$
ghost nodes (where $n_{g}$ depends on the support of the interpolation
and discretization stencils) that overlaps with adjacent blocks. During
the synchronization of the ghost nodes, the processes first prepare
the required data, {\it i.e.}, apply restriction and prediction operators
if required, and copy them into a buffer, which is then sent to the
receiving process. The receiver then extracts the data to the ghost
node layers of each block. We skip more details on this technical part in the interest of brevity. The handling of block-block interfaces {\it via} the ghost nodes is linked to the precise definition of blocks and grid, and it is critical for numerical stability. We elaborate on our choices in appendix~B (suppl. material)

Parallelization of an adaptive code is inherently more complicated
than for a code using static grids. The number of blocks can vary significantly
during a simulation and choosing the same number of processes throughout
the simulation is often not an efficient choice. On most supercomputers,
a simulation is interrupted after a given time limit is exceeded, which
is usually of the order of 24h. We exploit this fact to decide upon
resubmission of a simulation if the number of CPUs is increased, decreased
or can remain unchanged. We will discuss the performance in greater
detail in section~\ref{sec:Performance}.

\section{Numerical results and performance\label{sec:results}}

All simulations presented in this paper can be reproduced using the code\footnote{{https://github.com/adaptive-cfd/WABBIT}} and the parameter files from the supplementary material (See section E for an overview).

\subsection{Three vortices}
\begin{figure}

\includegraphics[width=0.9\textwidth]{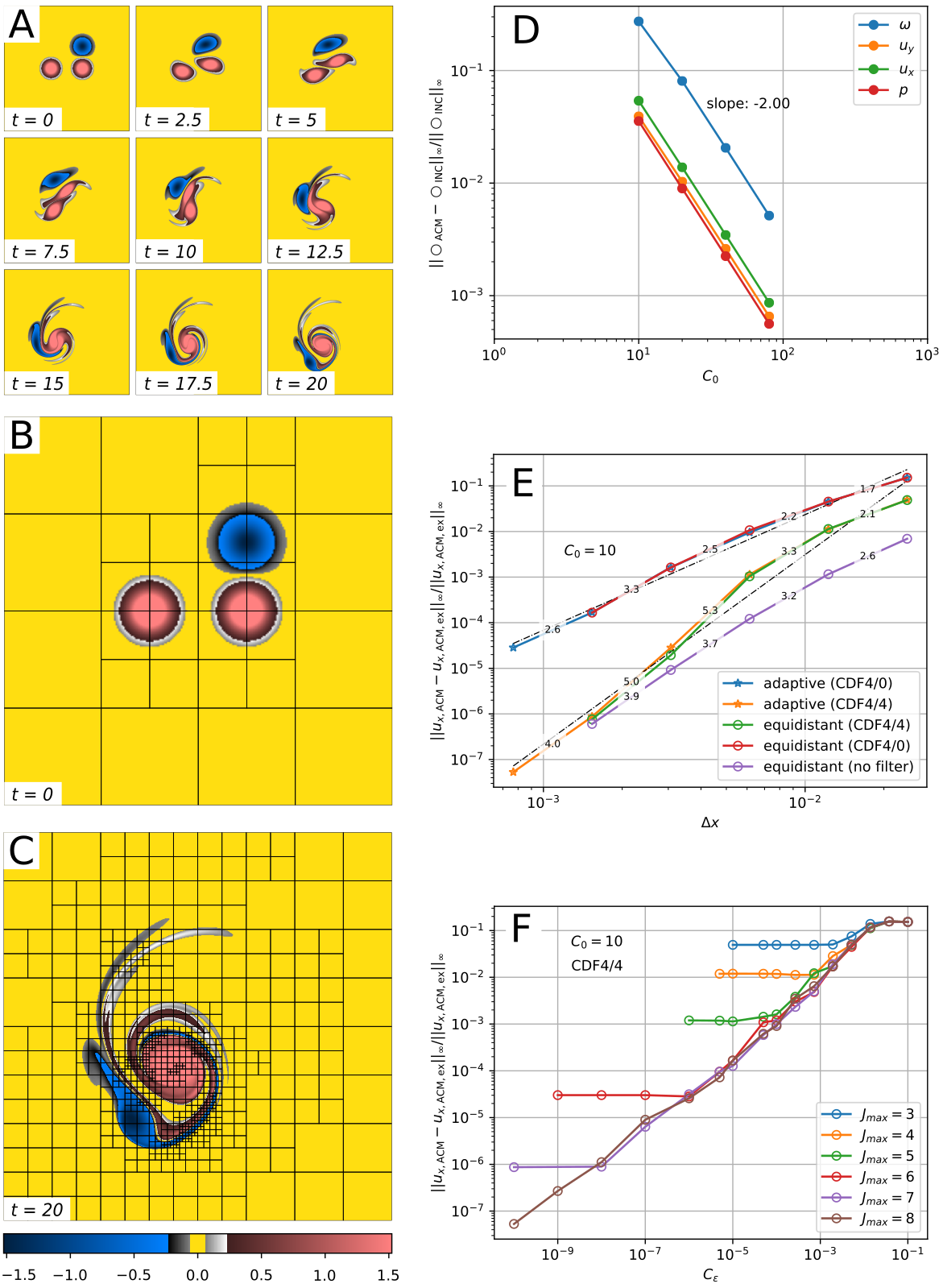}

\caption{Three vortices simulations. A-C: Time evolution of the vorticity field
in a simulation with $B_{s}=33$, $C_{\varepsilon}=10^{-5}$, $J_{max}=8$
and $C_{0}=10$ using the CDF 4/4 wavelet. Time evolution (A), initial condition (B) and final
time used for error calculation (C). D: Decay of the compressibility error with respect to the
incompressible solution as a function of $C_{0}$, computed using
a non-adaptive spectral discretization. Second order is observed.
E: Decay of the error w.r.t. the spectral ACM solution as a function
of $\Delta x$ for $C_{0}=10$. Shown are equidistant computations
with and without wavelet filtering and adaptive simulations with $C_{\varepsilon}$
small enough for the discretization error to be dominant. The
adaptive computations preserve the accuracy of the equidistant ones.
F: Decay of the error using CDF 4/4 w.r.t. the spectral ACM solution as a function
of $C_{\varepsilon}$ for different $J_{\mathrm{max}}$ and $C_{0}=10$.
Leftmost points of each line are also shown in E. \label{fig:Three-vortices-simulations.}}
\end{figure}

As a quantitative validation test for the fluid part in two-dimensions, we consider the vortex merger proposed in \cite{Kevlahan1997}. The setup involves no walls and is a good starting point for the validation of our new numerical method. At time $t=0$, three Gaussian vortices (Fig. \ref{fig:Three-vortices-simulations.}A,B) are set into the periodic domain of size $L=2\pi$, two with positive $(\Gamma=1$) and one with negative ($\Gamma=-1/2$) circulation. 
The initial vorticity field is given by
\[
\omega(x,t=0)\, = \, \sum_{i=1}^{3} \,\frac{\Gamma_{i}}{\pi\cdot a^{2}}\exp\left(-((x-x_{0,i})^{2}+(y-y_{0,i})^{2})/a^{2}\right),
\]
where $a=1/\pi$ and the vortex centers $(x_{0,i},y_{0,i})$ are located
at $(0.75,1)\pi$, $(1.25,1)\pi$ and $(1.25,1.25)\pi$. Note that
initially, the mean vorticity is small but not zero. The vortices
then interact nonlinearily, but their evolution remains deterministic.
Due to the low viscosity of $\nu=5\cdot10^{-5}$, thin vortex filaments
develop during the merging of the two positive vortices (Fig. \ref{fig:Three-vortices-simulations.}A).
We stop the simulation at $t=20$ (Fig. \ref{fig:Three-vortices-simulations.}C),
and use this time instant for comparing different simulations. The
creation of small scales due to the nonlinearity is an important test
for our adaptive framework.

We perform four types of numerical simulations for this problem:
\begin{enumerate}
\item Incompressible Navier--Stokes (INC) solved with a Fourier pseudospectral
method \cite{Engels2015a}, which is the quasi-exact reference solution,
\item ACM solved with the same Fourier pseudospectral method, as quasi-exact
solution of the artificial compressibility equations,
\item ACM solved with the fourth-order finite difference method on equidistant grids, with and without filtering at the finest scale
\item ACM solved with the fourth-order finite difference method on dynamically adaptive grids using multiresolution analysis with CDF 4/0 and CDF 4/4 wavelets. 
\end{enumerate}
Those simulations allow us to assess (i) the model error of the ACM,
(ii) the discretization error and (iii) the thresholding error. Combined,
they allow us to confirm the relation between the three numerical parameters
of the problem ($C_{o}$, $\Delta x$, $C_{\varepsilon}$) presented
in eq. (\ref{eq:scaling-no-penalization}). Instead of $\Delta x$,
it is however more practical to use the maximum refinement $J_{\mathrm{max}}$,
which is related to $\Delta x$ via $\Delta x=2^{-J_{\mathrm{max}}}L/(B_{s}-1)$.
In all ACM simulations of this problem, we set the pressure damping $C_{\gamma}=1$, as the sponge technique cannot be used for periodic domains.

\paragraph{(i) Compressibility error }

We first consider the Fourier pseudo-spectral solution of the artificial compressibility
eqns. (\ref{eq:u_eqn-1}-\ref{eq:p_eqn-1}). Comparing the spectral
ACM simulations with the incompressible reference solution from \cite{Kevlahan1997}
allows us to evaluate the compressibility error alone. 
We compute the reference solution with spatial resolution
of $1536\times1536$ modes using the \texttt{FLUSI} code
\cite{Engels2015a}, which uses the same Fourier pseudospectral
scheme as in \cite{Kevlahan1997}. Time integration is done with the classical
4th order Runge--Kutta scheme and $\mathrm{CFL}=0.15$. For all times
$t\geq0$, we verified that the isotropic Fourier spectrum decays
exponentially for large wavenumber $k$ below machine precision, hence both the reference solution and the spectral ACM simulations can be considered highly accurate.
Figure \ref{fig:Three-vortices-simulations.}D
shows the resulting convergence in $C_{0}$, confirming the second
order convergence from eq. (\ref{eq:ACM-model-error}).

\paragraph{(ii) Discretization error }

We next evaluate the discretization error of the finite difference
scheme on equidistant grids, fixing $C_{0}=10$ for computational
convenience, and compare with the corresponding spectral solution of the ACM equations.
The equidistant computations with resolutions from $N=256$ to $4096$
are performed using the \texttt{WABBIT} code with adaptivity turned
off. In this test case, the viscosity is sufficiently elevated to
stabilize the simulation without filtering at the maximum level $J_{\mathrm{max}}$.
Therefore we perform three types of equidistant computations shown
in Fig. \ref{fig:Three-vortices-simulations.}E, CDF 4/0 (red line),
CDF 4/4 (green) and without (purple) a filter. Without
filtering, the method is clearly fourth order accurate. If wavelet filtering at $J_\mathrm{max}$ is used, which corresponds to removing the detail coefficients, the error in the CDF 4/4 solution increases slightly, as expected, compared to the case without filtering.
Note that filtering is usually required for stabilization.
The convergence rate
locally varies, but a linear least squares fit to the data yields a slope of $4.04$.
With the CDF 4/4 wavelet, we can hence preserve the accuracy. The CDF 4/0
method on the other hand significantly increases the error and
deteriorates the convergence rate. We anticipate that this is due to the lack of scale separation in the CDF 4/0 wavelet. 
\paragraph{(iii) Thresholding error }

We now consider fully adaptive simulations, where we fix $C_{0}=10$
and use the refine -- evolve -- coarsen cycle in each time step
to remove aliasing errors. The block size $B_{s}$ is mainly a parameter
that impacts the computational efficiency, which is crucial only for
3D simulations. We here fix $B_{s}=33$ for all 2D simulations
and discuss the performance as a function of $B_{s}$ later in section
\ref{sec:Performance}. Two snapshots of an adaptive computation are shown
in Fig. \ref{fig:Three-vortices-simulations.}B,C. Visibly, the initial
grid contains less blocks than the terminal one. The refinement is
concentrated, as expected, in regions with sharp gradients.

In adaptive simulations, the thresholding error arises in addition
to discretization and compressibility errors. We compare with the
spectral ACM solution for $C_{0}$ fixed, hence the compressibility
error is absent. Figure~\ref{fig:Three-vortices-simulations.}F
shows the error as a function of $C_{\varepsilon}$ for different
$J_\mathrm{max}$. For large $C_{\varepsilon}$, all simulations run on the
coarsest possible level ($J_{\mathrm{max}}=1$) and therefore yield
the same error. When $C_{\varepsilon}$ is decreased, the error generally
decreases, until a saturation at a level-dependent value of $C_{\varepsilon,\mathrm{opt}}(J_{\mathrm{max}})$
is observed. At this point, the discretization- and thresholding error
are of the same size. Decreasing $C_{\varepsilon}$ further can thus
not increase the precision, as the discretization error then dominates.
Before the saturation, the error decays $\propto C_{\varepsilon}$.
Fig. \ref{fig:Three-vortices-simulations.}E shows the convergence
as a function of $\Delta x$ when using CDF 4/0 (blue line) and CDF 4/4 (orange), where we use the smallest
$C_{\varepsilon}$ used for each $J_{\mathrm{max}}$. The obtained
curves overlap the equidistant one almost perfectly -- we can hence
conclude that we preserve the accuracy and the order of the discretization scheme on adaptive grids. Consequently, all further results are obtained with the CDF 4/4 wavelet.

\subsection{A flapping wing}

\begin{figure}[h]
\begin{centering}
\includegraphics[width=0.85\textwidth]{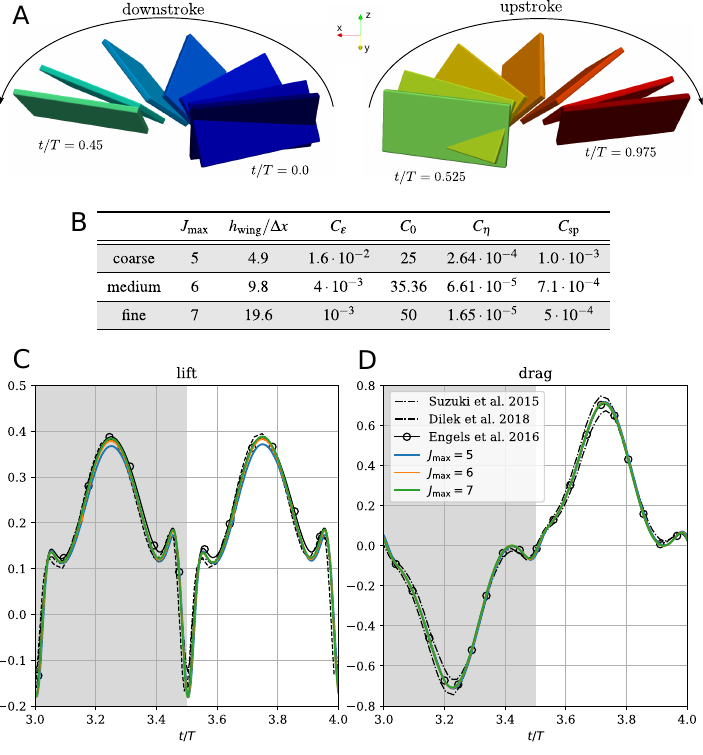}
\par\end{centering}
\caption{Suzuki's test case. \label{fig:Suzuki-test-case.}A: visualization of
the wingbeat. The rigid wing moves in an horizontal stroke plane with varying
angle of attack, the motion is symmetrical in down- and upstroke.
B: computational parameters used in the simulations. C-D:
Time evolution of lift and drag and comparison with results from the
literature. All results are obtained with CDF 4/4 wavelets.}
\end{figure}

As example of a time-dependent geometry in 3D, we consider the
test case proposed by Suzuki et al.~\cite{Suzuki2015}, Appendix B2.
It has also been considered as validation case in our previous work~\cite{Engels2015a} and by Dilek et al.~\cite{Dilek2019}. We therefore
limit the description of the setup to a minimum and refer the reader
to~\cite{Engels2015a} for details of the geometry and wingbeat kinematics.
The test consists of a single, rectangular, rigid flapping wing with finite
thickness $h_\mathrm{wing}=0.04171$ and length $R=1$. No inflow is imposed; the setup mimics hovering flight.
The wing moves in an horizontal stroke plane, its prescribed wingbeat
motion is visualized in Fig.~\ref{fig:Suzuki-test-case.}A. The two
half cycles, conventionally termed up- and downstroke, are symmetrical.
The size of the computational domain is chosen large enough to be
a good approximation for an unbounded flow, in our case we set $6R\times6R\times6R$,
which is slightly larger than~\cite{Suzuki2015,Dilek2019}. At the
outer border, we impose homogeneous Dirichlet conditions on $\underline{u}_{\infty}=0$
and $p_{\infty}=0$ using the sponge term in eqns.~(\ref{eq:u_eqn-1}-\ref{eq:p_eqn-1}).
The initial condition is $\underline{u}(t=0)=p(t=0)=0$.

As penalization is used for the wing, the scaling relations~(\ref{eq:scaling-penalized})
are used, and we use the constant $K_{\eta}=0.365$ justified from our previous
work \cite{Engels2015a}. Figure~\ref{fig:Suzuki-test-case.}B summarizes
the parameters used for a coarse, medium and fine resolution simulation.

The low Reynolds number of $Re=u_{\mathrm{tip}}c_{m}/\nu=100$,
which is comparable to that of a fruit fly, results in a smooth flow
topology. Here, $u_{\mathrm{tip}}$ is the mean wingtip velocity and
$c_{m}$ the chord length. The flow displays the characteristic features
of flapping flight, namely a strong leading edge vortex and a wingtip
vortex. The computational grid shows the expected refinement near
the wing, while the resolution in the remaining parts of the domain
is dictated by gradedness. Because of the absence of fine-scale structures
in the flow, the grid is in practice determined uniquely by the mask
function $\chi$, and the parameter $C_{\varepsilon}$ is not important.
Nonetheless, $C_{\varepsilon}$ has been scaled according to eq. (\ref{eq:scaling-penalized})
for consistency. We set a block size of $B_{s}=23$ (see section~\ref{sec:Performance}) and because $N_b$ does not strongly vary over time, the number of CPU is constant throughout the simulation.

Figure~\ref{fig:Suzuki-test-case.}C-D compares the lift- and drag force
obtained by present simulations with the results obtained by Suzuki
et al.~\cite{Suzuki2015}, Dilek et al \cite{Dilek2019} and our own previous
work~\cite{Engels2015a}. As already noted by Suzuki in their initial
publication using two different solvers, some discrepancy between
the reference solutions is observed. The difference, integrated over
the stroke period, is of the order of 5\%, which is why we do not
present convergence plots.

\subsection{Bumblebee in forward flight\label{sec:bumblebee}}

\begin{figure}
\begin{centering}
\includegraphics[width=0.85\textwidth]{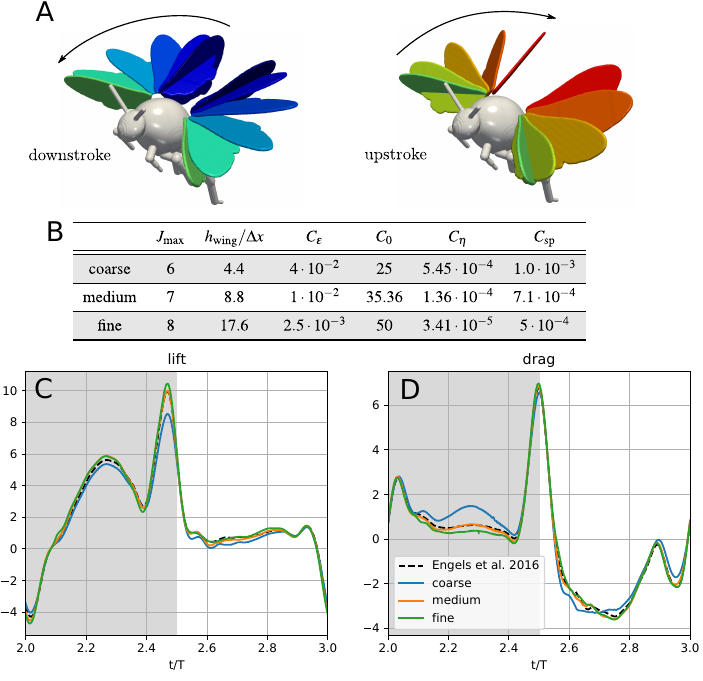}
\par\end{centering}
\caption{Bumblebee simulations. A: Wingbeat of the bumblebee model in forward
flight. During the downstroke, the angle of attack is more elevated
than during the upstroke. The stroke plane is inclined to generate
a forward force. B: computational parameters for the three levels of simulation. C-D: Lift- (=vertical) and drag (=horizontal) force. All results are obtained with CDF 4/4 wavelets.\label{fig:Bumblebee-1}}
\end{figure}

The last example in our hierarchy of validation cases is the simulation
of a model bumblebee. This model has also been used in our previous
work~\cite{Engels2019,Engels2015} and it was the first to be
studied in fully developed turbulence. In the interest of brevity, we limit the description of the model to a minimum and refer to ~\cite{Engels2019,Engels2015} for further details. The electronic supplementary material
to this article contains the parameter files used in the simulation for reproducibility. The penalization method allows for the inclusion of obstacles of arbitrary complexity, and we hence
include the bumblebees body. The wingbeat is illustrated in Fig.~\ref{fig:Bumblebee-1}A.
We are considering forward flight, which is why the stroke plane is
inclined with respect to horizontal. The wing kinematics differs significantly
between down- and upstroke. During the downstroke, also called power-stroke, most of the aerodynamic forces are generated, as the wing moves against the flow direction at elevated angle of
attack. The upstroke at much smaller angle of attack primarily serves
to enable the next downstroke. The wings are rigid and have a thickness $h_\mathrm{wing}=0.025R$. Our bumblebee has a Reynolds number
of $Re=u_{\mathrm{tip}}c_{m}/\nu=2000$ and as such produces a more
complicated flow topology than the previous test~\cite{Engels2015,Engels2017,Engels2019}.

The domain size is $L=8$, the block size is $B_{s}=23$ (see section
\ref{sec:Performance}) and the other numerical parameters of three simulations
levels are given in Fig.~\ref{fig:Bumblebee-1}B. All results are obtained with the CDF 4/4 wavelets.

\begin{figure}
\begin{centering}
\includegraphics[width=0.85\textwidth]{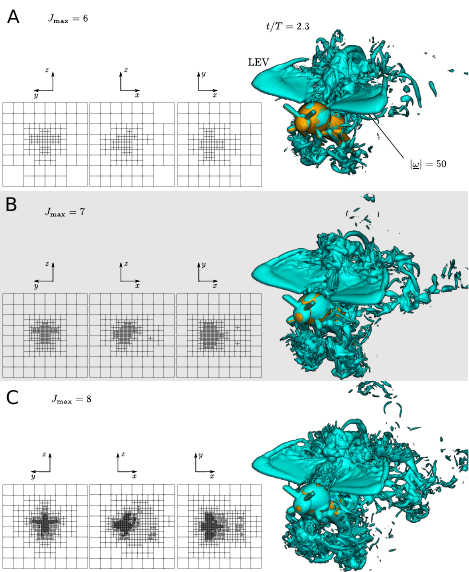}
\par\end{centering}
\caption{Bumblebee simulations. A (coarse), B (medium), C (fine resolution): Visualization of the flow field ($||\underline{\omega}||=50$ isosurface of vorticity magnitude) at $t/T=2.3$ (right) and 2D projections of the computational grid used at this time (left) for the three
simulations. All results are obtained with CDF 4/4 wavelets.\label{fig:Bumblebee-2}}
\end{figure}

The drag (Fig.~\ref{fig:Bumblebee-1}C) and lift (Fig.~\ref{fig:Bumblebee-1}D)
show good convergence and display features typical for flapping wings. The leading
edge vortex (LEV, cf. Fig. \ref{fig:Bumblebee-2}A) is built up during the downstroke and reaches its maximum intensity during the mid-downstroke ($t=2.25$), which results in a peak in the lift force. At the turning point ($t=2.5$) between down- and upstroke, the rapid wing rotation causes large peaks in
the forces.

Figure~\ref{fig:Bumblebee-2}A-C shows the same snapshot of the flow
field at $t/T=2.3$ in all three simulations. Shown is the $\left|\underline{\omega}\right|=50$
isosurface of vorticity magnitude. Visually, the three levels yield similar results for the large scale structures of the flow field. The grids concentrate, as expected, high resolution near the insect, but contrarily to the previous testcase parts of the wake are now likewise are refined to $J_{\mathrm{max}}$, because they are detected by the wavelet thresholding. In the coarse simulation, 3\% of the full grid is used, which decreases to 1.3\% and 0.6\% for the medium and fine resolution, respectively. 

\subsection{Performance\label{sec:Performance}}

For a specific simulation, three factors are important for the computational performance: the block size $B_s$, the number of blocks per CPU $N_b/N_{CPU}$ and the number of CPU itself. $B_s$ balances grid sparsity and data regularity and thus influences the solution, while the other parameters are purely computational. We clarify the influence of these parameters in this section.

The impact of $B_s$ can be divided into to contributions: the grid sparsity that depends on the physics of a simulation and the computational efficiency resulting from the data regularity. The latter depends only on the machine. We therefore study this part as follows. For each value of $B_s$, a random grid of given $N_b$ is first generated, then, 15 time steps of the ACM equations are performed on this grid. Each step consists of the refine-evolve-coarsening cycle explained in section \ref{subsec:Wavelet-Based-Multiresolution-Alg}. As we use random data, the adaptation is set to coarsen everywhere instead of using wavelet thresholding. After measuring the performance on this grid, a new, slightly denser grid is generated. The procedure is repeated until the allocated memory is exhausted.

\begin{figure}
\begin{centering}
\includegraphics[width=0.85\textwidth]{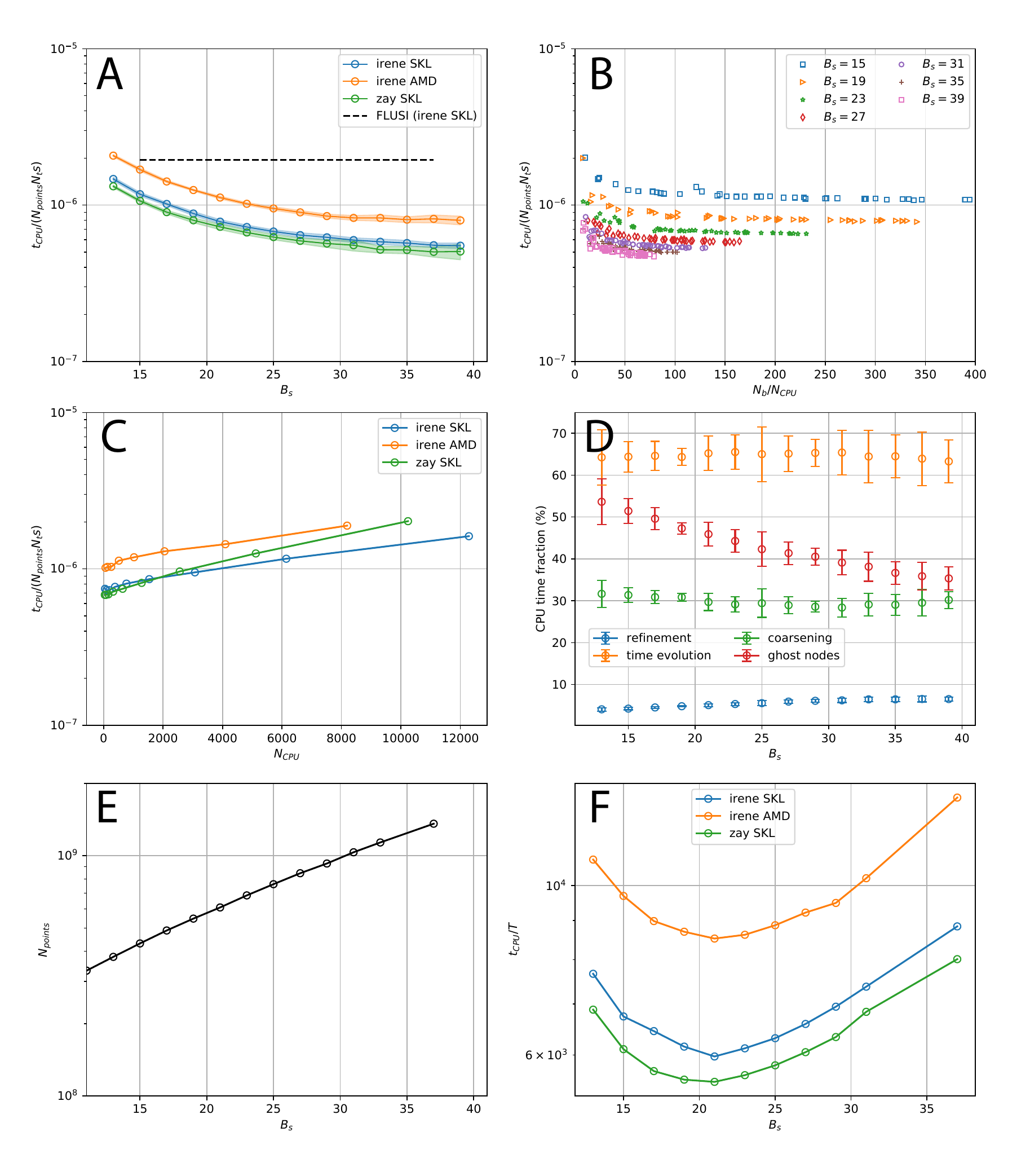}
\par\end{centering}
\caption{Computational performance. A-D: performance on randomly generated grids. E-F: performance of bumblebee simulations. A: CPU cost per grid point per right hand side evaluation as a function of $B_s$, measured on three different supercomputers (\texttt{Irene SKL}, \texttt{Irene AMD}, \texttt{Zay SKL}). For comparison, the cost of the \texttt{FLUSI} code on \texttt{Irene SKL} is included as well. Note only one simulation is shown, which does not depend on $B_{s}$. Shaded areas correspond to mean $\pm$ 1 s.d. 
B: CPU cost as a function of $N_b/N_{CPU}$ for different $B_s$, computed on $N_{CPU}=200$ cores on \texttt{Zay SKL}.
C: for $B_s=23$, CPU cost as a function of the number of CPU (weak scaling). For each datapoint $N_b/N_{CPU}$ is sufficiently high to be independent of that parameter.
D: fraction of time spend on grid refinement, time evolution and coarsening
(=100\%), as well as ghost node synchronization via \texttt{MPI}. Computed on $N_{CPU}=200$ cores on \texttt{Zay SKL}.
E: Bumblebee simulations, time-averaged total number of points on the grid $N_{\mathrm{points}}=N_{b}B_{s}^{3}$ as a function of $B_s$
F: Bumblebee simulations, total cost per time unit of a simulation with different $B_s$
All results are obtained with CDF 4/4 wavelets.\label{fig:performance}}
\end{figure}

Figure \ref{fig:performance}A shows the CPU time (in seconds) consumed per grid point and per right hand side evaluation, as the RK4 scheme has $s=4$ stages. This number is hence independent of the chosen time stepper, as our adaptivity only concerns spatial discretization. We used three different supercomputers, \texttt{Irene SKL}, \texttt{Irene AMD}\footnote{http://www-hpc.cea.fr/fr/complexe/tgcc.htm}, \texttt{Zay SKL}\footnote{http://www.idris.fr/}. For smaller $B_{s}$, the cost increases
because the overhead due to adaptivity becomes larger. Larger $B_{s}$
reduces the cost per point because of the increased data regularity.
It is important to compare the performance to other existing codes,
therefore the computational cost for the spectral code \texttt{FLUSI}~\cite{Engels2015a}
is also shown. We choose this code because we can simulate
exactly the same bumblebee setup with a completely different approach.
In this simulation, the resolution is $2048^{3}$ grid points and $N_{\mathrm{CPU}}=1152$.
The comparably low $N_{\mathrm{CPU}}$ ensures best possible performance
for the \texttt{FLUSI} code. The performance of \texttt{WABBIT} is about four times better than that of \texttt{FLUSI} for larger $B_{s}$, which demonstrates that the additional overhead
due to adaptivity is successfully mitigated by the block based data structure. 
Because of the full grid in the \texttt{FLUSI} case, the total
cost is much more elevated, but on the other hand it relies on highly optimized \texttt{FFT}
libraries. The similar performance per grid point and time step is thus a strong indicator 
for competitive performance of the \texttt{WABBIT} code.

Using dynamical grids the work load per CPU is not constant. It is therefore important to verify that
the performance remains indeed constant over a large range of blocks
per CPU, $N_{b}/N_{\mathrm{CPU}}$, as shown in Fig. \ref{fig:performance}B, computed using $N_{CPU}=200$ on \texttt{Zay SKL}. For $N_{b}/N_{\mathrm{CPU}}>50$,
the performance is nearly constant, while it drops for smaller values.

With increasing $N_{CPU}$, reasonable values for $N_{b}/N_{\mathrm{CPU}}$ are obtained only for very large simulations. At this point, the overhead of the global topology management can eventually become important and increases the cost per grid point (Fig. \ref{fig:performance}C). Improving this scaling behaviour is left for future work.

Figure~\ref{fig:performance}D shows the fraction of CPU time
spent on refinement, evolution and coarsening, which sum to 100\%.
The fractions spent on refinement, evolution and coarsening remain approximately constant over the range of $B_s$, while the ghost node synchronization and hence the parallelization overhead are much larger at small $B_s$. Synchronization is performed in all three of refinement, evolution and coarsening. Approximately 60\% --70\% are spent on computing the actual right hand side. We should note at this point that the ACM equations considered here are extremely simple and can be solved very efficiently. For more complicated equations, such as the full compressible Navier--Stokes equation, the portion of time spent on time evolution easily exceeds 90\%.

For the simulations presented earlier, we have not yet discussed how the block size $B_{s}=23$ has been chosen. Prior to performing the actual simulations, we performed a
series of runs with varying $B_{s}$ on the bumblebee configuration.
As the smallest $\Delta x$ depends on $B_{s}$, we cannot simply change $B_{s}$ without altering the
results. Therefore, we took advantage of the fact that for an exterior
flow, the domain size $L$ is an artifact from the simulation rather
than a physical quantity ($L\rightarrow\infty$ in reality). Thus,
for every $B_{s}$ we modified $L$ to obtain the same $\Delta x$.
All simulations run until $t=0.15$, and the other parameters are $J_{\mathrm{max}}=8$, $C_{\varepsilon}=10^{-3}$, $C_{0}=20$ and $C_{\eta}=1.65\cdot10^{-4}$.

Fig.~\ref{fig:performance}E
shows the total number of points on the grid, averaged over the time
span of the simulations. It is observed that as the computational
cost per point decreases (Fig.~\ref{fig:performance}A), the total number of points increases. Both
influences combined yield Fig.~\ref{fig:performance}F, which
shows the total cost per simulation time unit. It presents a minimum
for an intermediate value of $B_{s}$, where a compromise between
sparsity and data regularity is achieved. We therefore choose $B_{s}=23$,
but we also note that any value between $15$ and $30$ would have resulted in near-optimal performance

\section{Outlook: bumblebee in the wake of a fractal tree\label{sec:outlook}}

\begin{figure}
\begin{centering}
\includegraphics[width=0.85\textwidth]{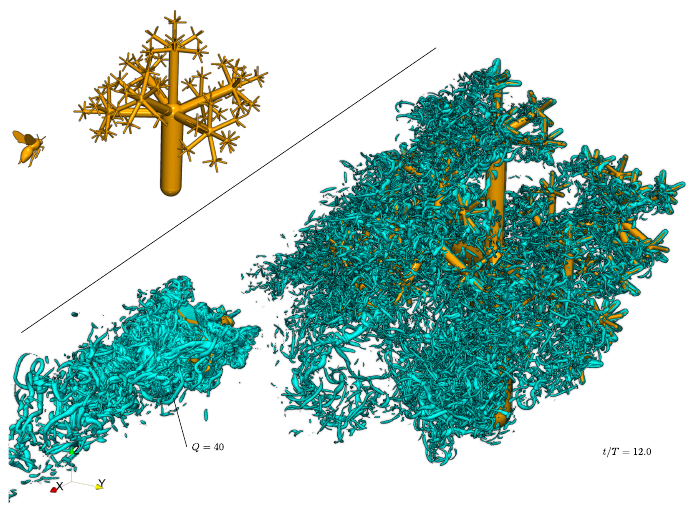}
\par\end{centering}
\caption{
Bumblebee behind a fractal tree. Top part shows the setup consisting of a bumblebee and a tree-inspired fractal turbulence generator composed of rigid cylinders. Bottom part illustrates the flow field with an isosurface of the $Q$-criterion.
\label{fig:bb-ft}
}
\end{figure}

To give an outlook and to illustrate the flexibility of the present adaptive wavelet approach to compute flows with multiple scales, we simulate the flow of the bumblebee model, described previously, in the wake of a model flower, corresponding to a bio-inspired turbulence generator. The latter is composed of an array of rigid cylinders arranged in a fractal manner, largely motivated by plants.  
For a detailed description of the geometry we refer the reader to \cite{Dreissigacker2017}. Both bumblebee and fractal tree are simulated using the volume penalization method in a large cubic domain of size $L=64R$, where $R$ is the wing length of the insect.
Figure \ref{fig:bb-ft} (top, left) illustrates the setup. 
The thinnest cylinders in the fractal tree have a diameter of $d_\mathrm{min}=0.114R$ and thus we have a Reynolds number of $\mathrm{Re}=u_\infty d_\mathrm{min} / \nu =250$. The other numerical parameters are chosen as in the coarse simulations of the bumblebee alone, described in section \ref{sec:bumblebee}.
The simulation is performed using CDF 4/4 wavelets where $C_\varepsilon=4 \cdot 10^{-2}$.
The flow field is visualized in Fig. \ref{fig:bb-ft} with an isosurface of the $Q$-criterion. It shows the turbulent wake produced by the tree model and the wake of the bumblebee after twelve wingbeats. The maximum number of refinement levels is $J_\mathrm{max}=9$ and the grid is composed, on average, of $N_b=2\cdot 10^{5}$ blocks with a total number of $2.5\cdot10^{9}$ grid points. The dynamically adapted grid allows us to use a large domain size, which would result in $11776^3$ points on the corresponding, uniform grid with the same $\Delta x$. We hence use on average only 0.15\% of the uniform grid. Note however that the domain size could possibly be reduced.
A detailed study on the flow physics of the bumblebee-plant interaction will be reported elsewhere.

\section{Conclusions\label{sec:conclusions}}

This article presents a novel framework for wavelet-based adaptive simulations of turbulent flows with a wide range of dynamically active scales in complex, time-dependent geometries on massively parallel computers, with a specific focus on flapping insect flight. 
It is implemented in the \texttt{WABBIT}-code, which is open-source in order to maximize its accessibility and its utility to the scientific community. 
We base our physical model on the original combination of artificial compressibility and volume penalization, and discretize the governing equations using classical centered finite differences on a locally uniform Cartesian grid together with a Runge--Kutta scheme for time integration, both of fourth order.
Motivated by the intermittent nature of our applications, we employ biorthogonal interpolating wavelets to introduce a dynamically adaptive grid tracking the solution in space and scale. In contrast to many previous wavelet-based implementations, the data structure is based on locally uniform blocks, which leads to an improved performance on modern CPUs. Thresholding the wavelet coefficients allows us to decide if a block can be coarsened at the end of a time step. As a block contains multiple wavelet coefficients but can only be coarsened as a whole, the number of points per block is a numerical parameter balancing the influence of data compression and CPU efficiency. The global grid topology is encoded using a custom tree-like data structure, as is natural in the wavelet-context. Even though our method has excellent compression properties, the applications we have in mind require many blocks and thus using massively parallel supercomputers is mandatory. To this end, blocks are distributed among MPI processes and the layer of overlapping ghost nodes is synchronized. This implementation yields high parallel efficiency using thousands of CPUs while maintaining both memory- and CPU time compression.

Through numerical estimations and a series of validation cases, we discussed how the parameters for this model have to be chosen in order to simultaneously reduce all occurring errors, {\it i.e.} compressibility, penalization, discretization and thresholding errors, while maintaining computational efficiency. We also assessed the impact of the block size on the performance and illustrated how it is chosen to yield maximum efficiency. Considering different supercomputers we furthermore showed that the performance per grid point and time step is very similar to that of the previously developed Fourier pseudospectral code \texttt{FLUSI}. Because this spectral code relies on the fast Fourier transform, a task for which highly tuned libraries are used, we can conclude that our new framework offers excellent performance and can indeed translate the sparsity of the grid into CPU time savings.

If we consider, for instance, the case of the high-resolution bumblebee computation with corresponding uniform grid resolution $5632^3$, we showed that the adaptive wavelet-based computation is, for similar accuracy, about 160 times faster than the classical Fourier pseudospectral code. While this is an optimistic scenario because the domain size can still be reduced, it still illustrates the strength of our novel approach.

We finally presented a simulation of an insect together with a fractal tree, which will be used in future work to study interactions of insects with their environment. Such multiscale simulations are out of reach for other approaches, in particular spectral codes, and we are not aware of other codes suitable for such problems. 

In future work we are planning to perform adaptive computations of insects using realistic bodies from micro CT data of real insects and mass-spring models for flexible wings \cite{Truong_etal_2020} to perform adaptive fluid-structure interaction simulations.

\section*{Acknowledgments}
TE, KS, MF gratefully acknowledge financial support from ANR (Agence Nationale de la Recherche), grant n\textsuperscript{\scriptsize o} 15-CE40-0019, and DFG (Deutsche Forschungsgemeinschaft), grant n\textsuperscript{\scriptsize o} SE 824/26-1, project AIFIT (Aerodynamics of Insect Flight In Turbulence). The authors were granted access to the HPC resources of IDRIS (Institut du d\'eveloppement et des ressources en informatique scientifique) under the allocation n\textsuperscript{\scriptsize o} 2018-91664 attributed by GENCI (Grand \'equipement national de calcul intensif). For this work we were also granted access to the HPC resources of Aix-Marseille Universit\'e financed by the project Equip{@}Meso (n\textsuperscript{\scriptsize o} ANR-10-EQPX- 29-01). The authors thankfully acknowledge financial support granted by MAEDI (Minist\`ere des Affaires \'etrang\`eres et du d\'eveloppement international), MEN (Minist\`ere de l'Education Nationale), MESRI (Minist\`ere de l'enseignement sup\'erieur, de la recherche et de l'innovation), DAAD (Deutscher Akademischer Austauschdienst) within the French-German program Procope, project FIFIT (Flapping Insect Flight In Turbulence). JR, TE gratefully acknowledge financial support from DFG, grant n\textsuperscript{\scriptsize o} SFB1029.

\bibliographystyle{plain}
\bibliography{bibliography}

\end{document}